\documentclass[11pt]{article}
\usepackage{amsfonts}
\usepackage{amsmath}
\usepackage{amssymb}
\usepackage{graphicx}
\usepackage[dvipsnames,usenames]{color}
\usepackage[pdfstartview=FitH ]{hyperref}
\usepackage[hyperpageref]{backref}
\begin{document}

\baselineskip 18pt
\def\o{\over}
\def\e{\varepsilon}
\title{\Large\bf  The\ \  Mean\ \  Square\ \  of\ \  Divisor\ \
 Function}
\author{Chaohua\ Jia\ and\ Ayyadurai\ Sankaranarayanan}
\date{}
\maketitle
{\small \noindent {\bf Abstract.} Let
$d(n)$ be the divisor function. In 1916, S. Ramanujan stated but
without proof that
$$
\sum_{n\leq x}d^2(n)=xP(\log x)+E(x),
$$
where $P(y)$ is a cubic polynomial in $y$ and
$$
E(x)=O(x^{{3\o 5}+\e}),
$$
where $\e$ is a sufficiently small positive constant. He also stated
that, assuming the Riemann Hypothesis(RH),
$$
E(x)=O(x^{{1\o 2}+\e}).
$$

In 1922, B. M. Wilson proved the above result unconditionally. The
direct application of the RH would produce
$$
E(x)=O(x^{1\o 2}(\log x)^5\log\log x).
$$
In 2003, K. Ramachandra and A. Sankaranarayanan proved the above
result without any assumption.

In this paper, we shall prove
$$
E(x)=O(x^{1\o 2}(\log x)^5).
$$
}

\vskip.1in
\noindent{\bf 1. Introduction}

Let $d(n)$ be the divisor function. In 1916, S. Ramanujan[9] stated
but without proof that
\begin{align*}
\qquad &\
d^2(1)+d^2(2)+d^2(3)+\cdots+d^2(n)\qquad\qquad\qquad\qquad\qquad
\qquad\ \ (1.1)\\
&=An(\log n)^3+Bn(\log n)^2+Cn\log n+Dn+ O(n^{{3\o 5}+\e}),
\end{align*}
here

\noindent--------------------------------

{\small 2010 Mathematics Subject Classification: Primary 11M;
Secondary 11M06.

}

$$
A={1\o \pi^2},\qquad\quad B={12\gamma-3\o \pi^2}-{36\o \pi^4}
\zeta'(2),
$$
where $\gamma$ is Euler's constant, $C,\,D$ are more complicated
constants, $\e$ is a sufficiently small positive constant. S.
Ramanujan[9] also stated that, assuming the Riemann Hypothesis(RH),
the error term in (1.1) can be improved to $O(n^{{1\o 2}+ \e})$.

Write
$$
E(x)=\sum_{n\leq x}d^2(n)-xP(\log x), \eqno (1.2)
$$
where
$$
P(x)=Ax^3+Bx^2+Cx+D.
$$
Then the statement of Ramanujan is that
$$
E(x)=O(x^{{3\o 5}+\e}), \eqno (1.3)
$$
and assuming the RH,
$$
E(x)=O(x^{{1\o 2}+\e}). \eqno (1.4)
$$

In 1922, B. M. Wilson[13] proved (1.4) unconditionally. By a general
theorem of M. K\"uhleitner and W. G. Nowak(see (5.4) in [5]), we
know
$$
E(x)=\Omega(x^{3\o 8}).  \eqno (1.5)
$$

Let $d_4(n)$ be the general divisor function which is the number of
representations of $n=d_1d_2d_3d_4$. In 1973, assuming
$$
\sum_{n\leq x}d_4(n)={1\o 6}x\log^3x+(2\gamma-{1\o 2})x\log^2x+ax
\log x+bx+O(x^\alpha),
$$
where $\gamma$ is Euler's constant, $a,\,b$ are constants, $\alpha$
is a constant strictly less than ${1\o 2}$, D. Suryanarayana and R.
Sitaramachandra Rao[10] proved
$$
E(x)=O(x^{1\o 2}\exp(-c(\log x)^{3\o 5}(\log\log x)^{-{1\o 5}})),
\eqno (1.6)
$$
where $c$ is a positive constant.

By Vinogradov's estimate, if ${T\o 2}\leq t\leq T$, then
$$
{1\o \zeta(1+2it)}\ll (\log T)^{2\o 3}(\log\log T)^{1\o 3}.
$$
So it is not difficult to prove
$$
E(x)=O(x^{1\o 2}(\log x)^{17\o 3}(\log\log x)^{1\o 3}). \eqno (1.7)
$$
The direct application of the RH (or even the quasi-RH) would produce
$$
E(x)=O(x^{1\o 2}(\log x)^5\log\log x). \eqno (1.8)
$$
In 2003, K. Ramachandra and A. Sankaranarayanan[8] proved (1.8)
without any assumption and put forward the following conjecture.

{\bf Conjecture}(Ramachandra-Sankaranarayanan). Assuming the RH, we
have
$$
E(x)=O(x^{1\o 2}).  \eqno (1.9)
$$

For the average situation, in 2005, H. Maier and A.
Sankaranarayanan[7] proved,
$$
{1\o X}\int_X^{2X}E^2(x)dx\ll X\exp(-c(\log X)^{3\o 5}(\log\log X)
^{-{1\o 5}}), \eqno (1.10)
$$
where $c$ is a positive constant.

In this paper, we shall prove the following theorem.

{\bf Theorem}. If $E(x)$ is defined in (1.2), then unconditionally
we have
$$
E(x)=O(x^{1\o 2}(\log x)^5). \eqno (1.11)
$$

Throughout this paper, we assume that $\e$ is a sufficiently small
positive constant and that $T$ is sufficiently large.

\vskip.3in
\noindent{\bf 2. Some lemmas}

{\bf Lemma 1}(Borel-Carath\'eodory). Suppose that $f(z)$ is
holomorphic in the disk $|z-z_0| \leq R$ and that in the circle
$z=z_0+Re^{i\theta}(0\leq\theta\leq 2\pi)$,
$$
{\rm Re}f(z)\leq M.
$$
Then in the disk $|z-z_0|\leq r(<R)$, we have
$$
|f(z)|\leq{2r\o R-r}\,M+{R+r\o R-r}\,|f(z_0)|.
$$

See Section 5.5 of [11].

{\bf Lemma 2}(Hadamard). Suppose that $f(z)$ is holomorphic in the
disk $|z-z_0|\leq R_3$, $R_1<R_2<R_3$. Write
$$
M_j=\max_{|z-z_0|=R_j}|f(z)|,\qquad\qquad\quad j=1,\,2,\,3.
$$
Then we have
$$
\log M_2\leq{\log({R_3\o R_2})\o \log({R_3\o R_1})}\cdot\log M_1
+{\log({R_2\o R_1})\o \log({R_3\o R_1})}\cdot\log M_3.
$$

See Section 5.3 of [11].

{\bf Lemma 3}. For $\alpha>0$ and $x>0$, we have
$$
{1\o 2\pi i}\int_{\alpha-i\infty}^{\alpha+i\infty}\Gamma(s)
x^{-s}ds=e^{-x}.
$$

See (2.15.2) in page 33 of [12].

{\bf Lemma 4}. For $-1\leq\sigma\leq 2$ and $|t|\geq 1$, we have
$$
\Gamma(\sigma+it)\ll|t|^{\sigma-{1\o 2}}e^{-{\pi\o 2}|t|}.
$$

See (4.12.2) in page 78 of [12].

{\bf Lemma 5}. For ${\rm Re}(s)>1$, let
$$
f(s)=\sum_{n=1}^\infty{a(n)\o n^s},
$$
where $a(n)=O(\psi(n))$, $\psi(n)$ is non-decreasing, and as
$\sigma\rightarrow 1^+$,
$$
\sum_{n=1}^\infty{|a(n)|\o n^\sigma}=O\Bigl({1\o (\sigma-1)^\alpha}
\Bigr).
$$
Then if $c>1$, $x$ is not an integer, and $N$ is the integer nearest
to $x$,
\begin{align*}
\sum_{n<x}a(n)&={1\o 2\pi i}\int_{c-iT}^{c+iT}f(s){x^s\o s}ds
+O\Bigl({x^c\o T(c-1)^\alpha}\Bigr)\\
&\quad +O\Bigl({\psi(2x)x\log x\o T}\Bigr)+O\Bigl({\psi(N)x\o
T|x-N|}\Bigr).
\end{align*}

See Lemma 3.12 in page 60 of [12].

{\bf Lemma 6}. For ${\rm Re}(s)>1$, we have
$$
\sum_{n=1}^\infty{d^2(n)\o n^s}={\zeta^4(s)\o \zeta(2s)}.
$$

See (1.2.10) in page 5 of [12].

{\bf Lemma 7}. For ${\rm Re}(s)\geq{1\o 2}$ and $|s-1|>1$, we have
$$
\zeta(s)=O \left ( |s| \right ).
$$

See (2.12.2) in page 29 of [12].

{\bf Lemma 8}. For $\sigma\geq 1$ and $t\geq 1$, we have
$$
{1\o \zeta(\sigma+it)}=O(\log t).
$$

See (3.11.8) in page 60 of [12].

{\bf Lemma 9}. For $t\geq 1$, we have
$$
\zeta({1\o 2}+it)=O(t^{{1\o 6}+\e}).
$$

See Theorem 5.5 in page 99 of [12].

{\bf Remark}. The bounds stated in Lemmas 8 and 9 suffice for our purpose
though better upper bounds are known.

{\bf Lemma 10}. For ${1\o 2}\leq\sigma\leq 1+\e$ and $t\geq 1$, we
have
$$
\zeta(\sigma+it)=O(t^{{1\o 3}(1-\sigma)+\e}).
$$

It follows from Lemma 9 and the explanation in Chapter 5 of [12].

{\bf Lemma 11}. We have
$$
\int_1^T|\zeta({1\o 2}+it)|^4dt=O(T\log^4T).
$$

See (7.6.1) in page 147 of [12].

{\bf Lemma 12}(Huxley). For $\sigma\geq{1\o 2}$, let $N(\sigma,\,T,
\,2T)$ denote the number of zeros $\rho=\beta+i\gamma$ of $\zeta(s)$
which satisfy $\beta\geq\sigma$ and $T\leq \gamma\leq 2T$. Then
$$
N(\sigma,\,T,\,2T)\ll T^{{12\o 5}(1-\sigma)+\e}.
$$

See [3].

{\bf Lemma 13}. For ${\rm Re}(z)>0$, we have
$$
\int_0^\infty e^{-zt}|\zeta({1\o 2}+it)|^2dt=2\pi e^{iz\o 2}\sum_
{l=1}^\infty d(l)\exp(2\pi ile^{iz})+f(z),
$$
where $f(z)$ is holomorphic in $|z|<4\e$.

This is Lemma 1 in [4].

Define
$$
D(s;\,{h\o k})=\sum_{l=1}^\infty{d(l)\o l^s}e(l{h\o k}). \eqno (2.1)
$$

{\bf Lemma 14}(Estermann). Suppose that $(h,\,k)=1$. The function
$D(s;\,{h\o k})$ is meromorphic in the whole plane with only one
pole of order 2 at $s=1$. In the neighborhood of $s=1$,
$$
D(s;\,{h\o k})={1\o k}\cdot{1\o (s-1)^2}+{2\o k}(\gamma-\log k)\cdot
{1\o (s-1)}+\cdots,
$$
where $\gamma$ is Euler's constant. At $s=0$, we have
$$
D(0;\,{h\o k})={1\o 4}-{1\o \pi
i}\sum_{a=1}^k\beta(a,\,k)\sum_{0<b<{k\o 2}}\eta(b,\,k)e(ab
{\overline{h}\o k}),
$$
where $\overline{h}h\equiv 1\,({\rm mod}\,k)$,
$$
\beta(a,\,k)=
\begin{cases} {1\o 1-e(-{a\o k})},\qquad\qquad &{\rm
if}\
1\leq a<k,\\
\ \ \ {1\o 2},&{\rm if}\ a=k,
\end{cases}
$$
and when $0<b<{k\o 2}$,
$$
0<\eta(b,\,k)<{1\o b}.
$$
Moreover, $D(s;\,{h\o k})$ satisfies the functional equation
$$
D(s;\,{h\o k})=2G^2(s)k^{1-2s}\Bigl(D(1-s;\,{\overline{h}\o k})-\cos
(\pi s)D(1-s;\,-{\overline{h}\o k})\Bigr),
$$
where
$$
G(s)=(2\pi)^{s-1}\Gamma(1-s).
$$

See (21), (34), (32), (29) and (19) in [1].

{\bf Lemma 15}. If $(m_1,\,m_2)=(n_1,\,n_2)=1$, then
$$
(m_1n_1^2,\,m_2n_2^2)=(m_1,\,n_2^2)(m_2,\,n_1^2).
$$

Proof. We have
$$
(m_1n_1^2,\,m_2n_2^2)=(m_1,\,n_2^2)\Bigl({m_1\o (m_1,\,n_2^2)}n_1^2,
\,m_2{n_2^2\o (m_1,\,n_2^2)}\Bigr).
$$
Since
$$
\Bigl({m_1\o (m_1,\,n_2^2)},\,m_2\Bigr)=1,\qquad\quad \Bigl({m_1\o
(m_1,\,n_2^2)},\,{n_2^2\o (m_1,\,n_2^2)}\Bigr)=1,
$$
we have
$$
\Bigl({m_1\o (m_1,\,n_2^2)}n_1^2,\,m_2{n_2^2\o (m_1,\,n_2^2)}\Bigr)
=\Bigl(n_1^2,\,m_2{n_2^2\o (m_1,\,n_2^2)}\Bigr)=(n_1^2,\,m_2).
$$
Thus, the conclusion of Lemma 15 follows.

{\bf Lemma 16}. If $a$ is a positive integer, then
$$
\sum_{M<m\leq 2^\e M}(m,\,a)\ll Md(a).
$$

Proof. We have
\begin{align*}
\sum_{M<m\leq 2^\e M}(m,\,a)&=\sum_{d|a}d\sum_{\substack{M< m\leq 2^\e M\\
(m,\,a)=d}} 1\\
&=\sum_{d|a}d\sum_{\substack{{M\o d}< m_1\leq{2^\e M\o d}\\
(m_1,\,{a\o d})=1}} 1\\
&\leq\sum_{d|a}d\sum_{{M\o d}<m_1\leq{2^\e M\o d}}1\\
&\ll\sum_{d|a}d\cdot{M\o d}=Md(a).
\end{align*}

{\bf Lemma 17}. Suppose that $0<A<B<2q$ and that $b$ is a positive
integer. Then
$$
\sum_{\substack{A<a\leq B\\ (a,\,q)=1\\ (a,\,b)=1}}
e(l{\overline{a}\o q}) \ll (l,\,q)^{1\o 2}q^{{1\o 2}+\e}b^\e.
$$
Here ${\overline {a}}$ is the integer such that $a{\overline {a}}
\equiv 1 ({\rm mod} \ q)$.

Proof. By Lemma 3 of [4], for $0<A<B<2q$, we have
$$
\sum_{\substack{A<a\leq B\\ (a,\,q)=1}} e(l{\overline{a}\o q})\ll
(l,\, q)^{1\o 2}q^{{1\o 2}+\e}.
$$
Hence,
\begin{align*}
\sum_{\substack{A<a\leq B\\ (a,\,q)=1\\ (a,\,b)=1}}
e(l{\overline{a}\o q}) &=\sum_{\substack{A<a\leq B\\
(a,\,q)=1}}\Bigl(\sum_{d|(a,\,b)}\mu(d)\Bigr)
e(l{\overline{a}\o q})\\
&=\sum_{d|b}\mu(d)\sum_{\substack{A<a\leq B\\ (a,\,q)=1\\
d|a}}e(l{\overline{a}\o q})\\
&=\sum_{d|b}\mu(d)\sum_{\substack{{A\o d}<t\leq{B\o d}\\
(dt,\,q)=1}} e(l\cdot{\overline{dt}\o q})\\
&=\sum_{\substack{d|b\\ (d,\,q)=1}}\mu(d)\sum_{\substack{{A\o d}<t\leq{B\o d}\\
(t,\,q)=1}} e(l\overline{d}\cdot{\overline{t}\o q})\\
&\ll\sum_{\substack{d|b\\
(d,\,q)=1}}|\mu(d)|\cdot(l\overline{d},\,q)^{1 \o 2}q^{{1\o 2}+\e}\\
&\ll(l,\,q)^{1\o 2}q^{{1\o 2}+\e}\sum_{d|b}1\\
&\ll(l,\,q)^{1\o 2}q^{{1\o 2}+\e}b^\e.
\end{align*}
Thus, Lemma 17 is proved.

\vskip.3in
\noindent{\bf 3. An asymptotic expression of $\bf
\zeta(1+it)$}

Let
$$
\rho_1=\beta_1+i\gamma_1,\ \ \rho_2=\beta_2+i\gamma_2,\ \ \cdots,\ \
\rho_J=\beta_J+i\gamma_J
$$
be all zeros of $\zeta(s)$ which satisfy $\beta\geq 1-4\e,\,T\leq
\gamma\leq 2T$. By Lemma 12,
$$
J=N(1-4\e,\,T,\,2T)\ll T^{11\e}. \eqno (3.1)
$$
We write domain $D$ as
$$
D=\{s=\sigma+it:\ \ 1-4\e\leq\sigma ,\ T\leq t\leq 2T\}.
$$
Write
\begin{align*}
&U_1=\bigcup_{j=1}^J(\gamma_j-(\log T)^{10},\
\gamma_j+(\log T)^
{10}),\\
&U_2=\bigcup_{j=1}^J(\gamma_j-2(\log T)^{10},\ \gamma_j+2(\log T)^
{10}),\\
&U_3=\bigcup_{j=1}^J(\gamma_j-3(\log T)^{10},\ \gamma_j+3(\log T)^
{10}),\\
\qquad\quad&U_4=\bigcup_{j=1}^J(\gamma_j-4(\log T)^{10},\
\gamma_j+4(\log T)^ {10}).\qquad\qquad\qquad\qquad  (3.2)
\end{align*}

After removing all domains of the form $\{s=\sigma+it:\ \ 1-4\e\leq
\sigma<1,\ t\in U_1\}$ in $D$, we denote the remained domain as
$D_1$. $D_1$ is a connected domain in which $\zeta(s)\ne 0$ so that
we can define a holomorphic function $\log \zeta(s)$ in $D_1$. For
${\rm Re}(s)>1$, Euler's product formula produces
$$
\log\zeta(s)=-\sum_p\log\Bigl(1-{1\o p^s}\Bigr)=\sum_p\sum_{m=1}^
\infty{1\o mp^{ms}}=\sum_{n=2}^\infty{\Lambda_1(n)\o n^s}, \eqno
(3.3)
$$
where
$$
\Lambda_1(n)={\Lambda(n)\o \log n}.
$$

After removing all domains of the form $\{s=\sigma+it:\ \
1-4\e\leq\sigma,\ t\in U_2\}$ in $D$, we denote the remained domain
as $D_2$. Now Lemma 1 can be applied. Take $f(z) =\log\zeta(z)$. For
$s=\sigma+it\in D_2,\,1-2\e\leq\sigma\leq 2$, let the center of
circle be $z_0=2+it$, the radius of bigger circle be $R=2-(1-4\e)=
1+4\e$, the radius of smaller circle be $r=2-(1-2\e)=1+2\e$. On the
bigger circle, by Lemma 7,
$$
{\rm Re}\log\zeta(z)=\log|\zeta(z)|\leq C\log T,
$$
where $C$ is a positive constant. Thus, for $s$ in the smaller circle,
Lemma 1 yields
$$
|\log\zeta(s)|\leq{2r\o R-r}\cdot C\log T+{R+r\o R-r}\cdot|\log
\zeta(2+it)|\ll \log T.
$$
For ${\rm Re}(s)\geq 2$, it is easy to see
$$
\log\zeta(s)=O(1).
$$
Hence, for $s=\sigma+it\in D_2,\,\sigma\geq 1-2\e$, we have
$$
|\log\zeta(s)|\ll\log T. \eqno (3.4)
$$

After removing all domains of the form $\{s=\sigma+it:\ \
1-4\e\leq\sigma,\ t\in U_3\}$ in $D$, then limiting $\sigma\geq
1-2\e$, we denote the obtained domain as $D_3$. Now Lemma 2 can be
applied. Take $f(z)=\log\zeta(z)$. For $s=\sigma+it\in D_3,\,1-
\e\leq\sigma\leq 1+\e$, let the center of circle be $z_0=2+it,\,R_3
=2-(1-2\e)=1+2\e,\, R_2=2-(1-\e)=1+\e,\,R_1=2-(1+\e)=1-\e$. By
(3.4), $M_3\ll\log T$. It is obvious that $M_1=O(1)$. Lemma 2 yields
\begin{align*}
\log M_2&\leq{\log({1+2\e\o 1+\e})\o \log({1+2\e\o 1-\e})}\cdot\log
M_1+{\log({1+\e\o 1-\e})\o \log({1+2\e\o 1-\e})}\cdot\log M_3\\
&\leq O(1)+{2\e+O(\e^2)\o 3\e+O(\e^2)}\cdot\log\log T\\
&=O(1)+({2\o 3}+O(\e))\log\log T\\
&\leq{3\o 4}\log\log T.
\end{align*}
Hence, for $s=\sigma+it\in D_3,\,1-\e\leq\sigma\leq 1+\e$, we have
$$
|\log\zeta(s)|\leq(\log T)^{3\o 4}.
$$
For ${\rm Re}(s)\geq 1+\e$, it is obvious that
$$
{1\o \zeta(s)}=O_{\e}(1).
$$
Thus, for $s=\sigma+it\in D_3,\,\sigma\geq1-\e$, we have
$$
{1\o \zeta(s)}\ll\exp((\log T)^{3\o 4}). \eqno (3.5)
$$

After removing all domains of the form $\{s=\sigma+it:\ \
1-4\e\leq\sigma,\ t\in U_4\}$ in $D$, then limiting $\sigma\geq
1-\e$, we denote the obtained domain as $D_4$. For $s\in D_4,\,u\geq
0,\,|v|\leq(\log T)^3$, we have
$$
{1\o \zeta(s+u+iv)}\ll\exp((\log T)^{3\o 4}). \eqno (3.6)
$$

For $s=1+it\in D_4,\,w=u+iv,\,X>1$, we have
\begin{align*}
&\ \,{1\o 2\pi i}\int_{u=\e,\,|v|\leq(\log T)^3}{1\o
\zeta(s+w)}\cdot\Gamma(w)X^wdw\\
&={1\o 2\pi i}\int_{u=\e,\,|v|\leq(\log
T)^3}\sum_{n=1}^\infty{\mu(n)
\o n^{s+w}}\cdot\Gamma(w)X^wdw\\
&=\sum_{n=1}^\infty{\mu(n)\o n^s}\cdot{1\o 2\pi i}\int_{u=\e,\,|v|
\leq(\log T)^3}\Gamma(w) \left ( {X\o n} \right )^wdw.
\end{align*}
By Lemma 4, if $|v|\geq 1$, then on the vertical line $u=\epsilon$, we have
$$
\Gamma(w)\ll |v|^{\e-{1\o 2}}e^{-{\pi\o 2}|v|}.
$$
Hence,
\begin{align*}
&\ \,{1\o 2\pi i}\int_{u=\e,\,|v|>(\log T)^3}\Gamma(w)({X\o n})^w
dw\\
&\ll({X\o n})^\e\int_{u=\e,\,|v|>(\log T)^3}|\Gamma(w)||dw|\\
&\ll({X\o n})^\e\int_{|v|>(\log T)^3}|v|^{\e-{1\o 2}}
e^{-{\pi\o 2}|v|}dv\\
&\ll({X\o n})^\e\int_{(\log T)^3}^\infty e^{-{\pi\o 2}v}dv\\
&\ll({X\o n})^\e\exp(-{\pi\o 2}(\log T)^3).
\end{align*}
By Lemma 3,
$$
{1\o 2\pi i}\int_{\e-i\infty}^{\e+i\infty}\Gamma(w)({X\o n})^wdw=
e^{-{n\o X}}.
$$
Therefore it follows that
\begin{align*}
&\ \,{1\o 2\pi i}\int_{u=\e,\,|v|\leq(\log T)^3}{1\o \zeta(s+w)}
\cdot\Gamma(w)X^wdw\\
&=\sum_{n=1}^\infty{\mu(n)\o n^s}\Bigl(e^{-{n\o X}}+O(({X\o n})^\e
\exp(-{\pi\o 2}(\log T)^3))\Bigr)\\
&=\sum_{n=1}^\infty{\mu(n)\o n^s}e^{-{n\o X}}+O(X^\e\exp(-{\pi\o
2}(\log T)^3)).
\end{align*}

We move the line of integration to ${\rm Re}(w)=-\e$. At $w=0$,
$\Gamma(w)$ has a pole of order 1 with residue 1. Hence, the residue
of ${1\o \zeta(s+w)}\cdot\Gamma(w)X^w$ at $w=0$ is ${1\o \zeta(s)}$.
In two horizontal lines , by (3.6),
\begin{align*}
&\ \,{1\o 2\pi i}\int_{-\e\leq u\leq\e,\,|v|=(\log T)^3}{1\o
\zeta(s+w)}\cdot\Gamma(w)X^wdw\\
&\ll X^\e\exp((\log T)^{3\o 4})\int_{-\e}^\e e^{-{\pi \o 2}(\log
T)^3}du\\
&\ll X^\e\exp(-(\log T)^3).
\end{align*}
The integration on ${\rm Re}(w)=-\e$ is
\begin{align*}
&\ \,{1\o 2\pi i}\int_{u=-\e,\,|v|\leq(\log T)^3}{1\o \zeta(s+w)}
\cdot\Gamma(w)X^wdw\\
&\ll X^{-\e}\exp((\log T)^{3\o 4})\Bigl(\int_{u=-\e,\,|v|\leq(\log
T)^3}|\Gamma(w)||dw|\Bigr)\\
&\ll X^{-\e}\exp((\log T)^{3\o 4})\Bigl(\int_{u=-\e,\,|v|\leq 1}
|\Gamma(w)||dw|\\
&\ +\int_{u=-\e,\,1\leq |v|\leq(\log T)^3}|\Gamma(w)||dw|\Bigr)\\
&\ll X^{-\e}\exp((\log T)^{3\o 4})\Bigl(\int_{u=-\e,\,|v|\leq 1}
{|dw|\o |w|} \\
&\ +\int_{1\leq |v|\leq(\log T)^3}|v|^{-\e-{1\o 2}}e^{-{\pi\o
2}|v|}dv\Bigr)\\
&\ll_\e X^{-\e}\exp((\log T)^{3\o 4}).
\end{align*}

Combining all of the above, we get (with $s=1+it$)
\begin{align*}
\qquad \qquad{1\o \zeta(s)}&=\sum_{n=1}^\infty{\mu(n)\o n^s}e^{-{n\o
X}}+O(X^\e
\exp(-(\log T)^3))\qquad\qquad\qquad (3.7)\\
&\qquad\qquad+O(X^{-\e}\exp((\log T)^{3\o 4})).
\end{align*}
Therefore we obtain an asymptotic expression of $\zeta(1+it)$ as
follows.

{\bf Proposition 1}. Suppose that $T\leq t\leq 2T,\,t\not\in U_4$
and
$$
X=\exp({2\o \e}(\log T)^{3\o 4}).  \eqno (3.8)
$$
Then we have
$$
{1\o \zeta(1+it)}=\sum_{n\leq X}{\mu(n)\o n^{1+it}}e^{-{n\o
X}}+O(1). \eqno (3.9)
$$

\vskip.3in
\noindent{\bf 4. A mean value estimate on $\bf \zeta(s)$}

In this section, we shall prove the following mean value estimate on
$\zeta(s)$.

{\bf Proposition 2}. If $k$ is any given positive number, then we have
$$
\int_1^T{|\zeta({1\o 2}+it)|^4\o |\zeta(1+2it)|^k}dt\ll_k T\log^4T.
$$

Firstly we shall prove the following Proposition 3. We use the
method of Iwaniec[4] essentially but with some modification and
refinement.

{\bf Proposition 3}. Suppose that $N\ll T^{{1\o 16}-\e}$ and that
for $N<n\leq 2^\e N$, $a(n)=O(N^{-1+\e})$. Then
$$
\int_{T\o 2}^T|\zeta({1\o 2}+it)|^4\Bigl|\sum_{N< n\leq 2^\e
N}{a(n)\o n^{2it}}\Bigr|^2dt\ll {T\log^4T\o N^{1-8\e}}.
$$

Proof. By the discussion in Section 2 of [4], we shall estimate
\begin{align*}
&\log T\sum_{r\leq{1\o 2\e\log 2}\log T+O(1)}\int_0^\infty
e^{-{t\o T}}|\zeta({1\o 2}+it)|^2\cdot\\
&\cdot\Bigl|\sum_{2^{\e r}< m\leq 2^\e\cdot 2^{\e r}}{1\o m^{{1\o
2}+it}}\Bigr|^2\Bigl|\sum_{N< n\leq 2^\e N}{a(n)\o
n^{2it}}\Bigr|^2dt.
\end{align*}
Write
\begin{align*}
&\ \,\Bigl|\Bigl(\sum_{M< m\leq 2^\e M}{1\o m^{{1\o 2}+it}}\Bigr)
\Bigl(\sum_{N< n\leq 2^\e N}{a(n)\o n^{2it}}\Bigr)\Bigr|^2\\
&=\Bigl|\sum_{K< k\leq 8^\e K}{b(k)\o k^{it}}\Bigr|^2=\sum_{K<k,\,h
\leq 8^\e K}b(k)\overline{b(h)}({h\o k})^{it},
\end{align*}
where $M=2^{\e r},\,M\ll T^{1\o 2},\,K=MN^2$,
$$
b(k)=\sum_{\substack{mn^2=k\\ M< m\leq 2^\e M\\ N< n\leq 2^\e
N}}{a(n) \o m^{1\o 2}}.
$$
In the following we shall estimate
\begin{align*}
&\ \,\int_0^\infty e^{-{t\o T}}|\zeta({1\o 2}+it)|^2\Bigl|\sum_{K<
k\leq 8^\e K}{b(k)\o k^{it}}\Bigr|^2 \qquad\qquad\qquad\qquad\qquad (4.1)\\
\qquad&=\sum_{K< k,\,h\leq 8^\e K}b(k)\overline{b(h)}\int_0^\infty
e^{-({1 \o T}-i\log({h\o k}))t}|\zeta({1\o 2}+it)|^2dt.
\end{align*}

Let
$$
z={1\o T}-i\log({h\o k}) \eqno (4.2)
$$
and note that
$$
|z| \leq  {1\o T} + \Bigl|\log(\frac {h}{k})\Bigr| < 4\e
$$
for $K < k,h \leq 8^\e K$.

By Lemma 13,
\begin{align*}
&\ \,\int_0^\infty e^{-zt}|\zeta({1\o 2}+it)|^2dt\\
\qquad &=2\pi e^{i\o 2T}({h\o k})^{1\o 2}\sum_{l=1}^\infty
d(l)\exp(2\pi
il({h\o k})e^{i\o T})+O(1)\qquad\qquad\qquad\qquad (4.3)\\
&=2\pi e^{i\o 2T}({h\o k})^{1\o 2}\sum_{l=1}^\infty d(l)e(l\cdot{h
\o k})\exp(2\pi il({h\o k})(e^{i\o T}-1))+O(1)\\
&=2\pi e^{i\o 2T}({h\o k})^{1\o 2}\sum_{l=1}^\infty d(l)e(l\cdot{h
\o k})\exp(2\pi ilx)+O(1),
\end{align*}
where
$$
x={h\o k}(e^{i\o T}-1). \eqno (4.4)
$$

The contribution of the term $O(1)$ to (4.1) is
\begin{align*}
&\ \,O\Bigl(\sum_{K< k,\,h\leq 8^\e K}|b(k)b(h)|\Bigr)\\
&\ll\sum_{M<m_1\leq 2^\e M}{1\o m_1^{1\o 2}}\sum_{M<m_2\leq 2^\e M}
{1\o m_2^{1\o 2}}\sum_{N<n_1\leq 2^\e N}|a(n_1)|\sum_{N<n_2\leq 2^
\e N}|a(n_2)|\\
&\ll M\sum_{N<n_1\leq 2^\e N}{1\o N^{1-\e}}\sum_{N<n_2\leq 2^\e N}
{1\o N^{1-\e}}\\
&\ll MN^{2\e}\ll {T\o N^{1-8\e}}.
\end{align*}

Let
$$
S(x;\,{h\o k})=\sum_{l=1}^\infty d(l)e(l\cdot{h\o k})\exp(2\pi ilx).
\eqno (4.5)
$$
Write
$$
\frak{z}=-2\pi ix=4\pi({h\o k})\sin({1\o 2T})e^{i\o 2T}. \eqno (4.6)
$$
By the discussion in Section 3 of [4], we know
$$
S(x;\,{h\o k})={1\o 2\pi i}\int_{1+\e-i\infty}^{1+\e+i\infty}D(s;\,
{h\o k})\Gamma(s)\frak{z}^{-s}ds, \eqno (4.7)
$$
where
$$
D(s;\,{h\o k})=\sum_{l=1}^\infty{d(l)\o l^s}e(l{h\o k}).
$$

In the following we write
$$
k^\ast={k\o (k,\,h)},\qquad\qquad h^\ast={h\o (k,\,h)}. \eqno (4.8)
$$
We move the line of integration from ${\rm Re}(s)=1+\e$ to ${\rm
Re}(s)=-\e$, and get
\begin{align*}
S(x;\,{h\o k})&={1\o 2\pi i}\int_{1+\e-i\infty}^{1+\e+i\infty}D(s;\,
{h^\ast\o k^\ast})\Gamma(s)\frak{z}^{-s}ds\\
&={1\o 2\pi i}\int_{-\e-i\infty}^{-\e+i\infty}D(s;\,{h^\ast\o
k^\ast})\Gamma(s)\frak{z}^{-s}ds+R_1(T;\,h,\,k)+R_0(T;\,h,\,k)\ (4.9)\\
&=R(T;\,h,\,k)+R_1(T;\,h,\,k)+R_0(T;\,h,\,k),
\end{align*}
where
$$
R(T;\,h,\,k)={1\o 2\pi
i}\int_{-\e-i\infty}^{-\e+i\infty}D(s;\,{h^\ast \o
k^\ast})\Gamma(s)\frak{z}^{-s}ds, \eqno (4.10)
$$
$R_1(T;\,h,\,k)$ and $R_0(T;\,h,\,k)$ are residues of
$D(s;\,{h^\ast\o k^\ast})\Gamma(s)\frak{z}^{-s}$ coming from the poles
at $s=1$ and $s=0$ respectively.

By the discussion in Section 3 of [4] and Lemma 14, we know that
\begin{align*}
\qquad\qquad R_1(T;\,h,\,k)&={1\o
\frak{z}k^\ast}(\gamma-\log\frak{z}-2\log k^
\ast)\ll {T\log T\o k^\ast},\qquad\qquad (4.11)\\
R_0(T;\,h,\,k)&=D(0;\,{h^\ast\o k^\ast})\qquad\qquad\qquad\qquad\qquad
\qquad\qquad\quad\ \ (4.12)\\
&={1\o 4}-{1\o \pi
i}\sum_{a=1}^{k^\ast}\beta(a,\,k^\ast)\sum_{0<b<{k^\ast\o
2}}\eta(b,\,k^\ast)e(ab {\overline{h^\ast}\o k^\ast}).
\end{align*}
Now we see the contribution of $R_1(T;\,h,\,k),\,R(T;\,h,\,k)$ and
$R_0(T;\,h,\,k)$ to (4.1).

1. The contribution of $R_1(T;\,h,\,k)$

We note that $\frac {h}{k} \ll 1$ for $K < h,k \le 8^{\e} K$. Therefore
the contribution of $R_1(T;\,h,\,k)$ is
\begin{align*}
&\ll\sum_{K<k,\,h\leq 8^\e K}|b(k)b(h)||R_1(T;\,h,\,k)|\\
&\ll\sum_{K<k,\,h\leq 8^\e K}|b(k)b(h)|\cdot{T\log T\o k}(k,\,h)
\end{align*}
\begin{align*}
&\ll T\log T\sum_{M< m_1\leq 2^\e M}\sum_{M< m_2\leq 2^\e M}\sum_{N
< n_1\leq 2^\e N}\sum_{N< n_2\leq 2^\e N}{|a(n_1)|\o m_1^{1\o 2}}
\cdot\\
\qquad &\qquad\qquad\cdot{|a(n_2)|\o m_2^{1\o 2}}\cdot{1\o m_1n_1^2}
(m_1n_1^2,\,m_2n_2^2)\qquad\qquad\qquad\qquad\qquad\ \ (4.13)\\
&\ll T\log T\cdot{1\o MN^{2-2\e}}\cdot{1\o MN^2}\sum_{M< m_1\leq 2^
\e M}\sum_{M< m_2\leq 2^\e M}\cdot\\
&\qquad\qquad\cdot\sum_{N< n_1\leq 2^\e N}\sum_{N< n_2\leq 2^\e N}
(m_1n_1^2,\,m_2n_2^2)\\
&={T\log T\o M^2N^{4-2\e}}\sum_{M< m_1\leq 2^\e M}\sum_{M< m_2\leq
2^\e M}\sum_{N< n_1\leq 2^\e N}\sum_{N< n_2\leq 2^\e N}(m_1n_1^2,
\,m_2n_2^2).
\end{align*}

By Lemmas 15 and 16,
\begin{align*}
&\ \,\sum_{M< m_1\leq 2^\e M}\sum_{M< m_2\leq 2^\e M}\sum_{N< n_1
\leq 2^\e N}\sum_{N< n_2\leq 2^\e N}(m_1n_1^2,\,m_2n_2^2)\\
&=\sum_{d\leq 2^\e M}\sum_{M< m_1\leq 2^\e M}\sum_{\substack{M<
m_2\leq 2^\e M\\ (m_1,\,m_2)=d}}\sum_{r\leq 2^\e N}\sum_{N< n_1\leq
2^\e N}\cdot\\
&\qquad\qquad\cdot\sum_{\substack{N< n_2\leq 2^\e N\\
(n_1,\,n_2)=r}}(m_1n_1^2,\,m_2n_2^2)\\
&=\sum_{d\leq 2^\e M}d\sum_{{M\o d}< m_1'\leq {2^\e M\o
d}}\sum_{\substack{{M\o d}< m_2'\leq {2^\e M\o d}\\
(m_1',\,m_2')=1}}\sum_{r\leq 2^\e N}r^2\sum_{{N\o r}< n_1'\leq
{2^\e N\o r}}\cdot\\
&\qquad\qquad\cdot\sum_{\substack{{N\o r}< n_2'\leq {2^\e N\o r}\\
(n_1',\,n_2')=1}}(m_1'n_1'^2,\,m_2'n_2'^2)\\
&=\sum_{d\leq 2^\e M}d\sum_{{M\o d}< m_1'\leq {2^\e M\o
d}}\sum_{\substack{{M\o d}< m_2'\leq {2^\e M\o d}\\
(m_1',\,m_2')=1}}\sum_{r\leq 2^\e N}r^2\sum_{{N\o r}< n_1'\leq {2^\e N\o r}}\cdot\\
&\qquad\qquad\cdot\sum_{\substack{{N\o r}< n_2'\leq {2^\e N\o r}\\
(n_1',\,n_2')=1}} (m_1',\,n_2'^2)(m_2',\,n_1'^2)\\
&\leq\sum_{d\leq 2^\e M}d\sum_{{M\o d}< m_1'\leq {2^\e M\o d}}\sum_
{{M\o d}< m_2'\leq {2^\e M\o d}}\sum_{r\leq 2^\e N}r^2\sum_{{N\o r}
< n_1'\leq {2^\e N\o r}}\cdot\\
&\qquad\qquad\cdot\sum_{{N\o r}<n_2'\leq {2^\e N\o
r}}(m_1',\,n_2'^2) (m_2',\,n_1'^2)
\end{align*}
\begin{align*}
&=\sum_{d\leq 2^\e M}d\sum_{r\leq 2^\e N}r^2\sum_{{N\o r}<n_2'\leq
{2^\e N\o r}}\sum_{{M\o d}< m_1'\leq {2^\e M\o
d}}(m_1',\,n_2'^2)\cdot\\
&\qquad\qquad\cdot\sum_{{N\o r}< n_1'\leq {2^\e N\o r}}\sum_{{M\o
d}< m_2'\leq {2^\e M\o d}}(m_2',\,n_1'^2)\\
&\ll\sum_{d\leq 2^\e M}d\sum_{r\leq 2^\e N}r^2\sum_{{N\o r}<n_2'\leq
{2^\e N\o r}}{M\o d}\cdot d(n_2'^2)\sum_{{N\o r}<n_1'\leq {2^\e N\o
r}}{M\o d}\cdot d(n_1'^2)\\
&\ll_{\e} M^2N^{2\e}\sum_{d\leq 2^\e M}{1\o d}\sum_{r\leq 2^\e N}r^2
\sum_{{N\o r}<n_1'\leq {2^\e N\o r}}\sum_{{N\o r}<n_2'\leq {2^\e N\o
r}}1\\
&\ll M^2N^{2\e}\log(2M)\sum_{r\leq 2^\e N}r^2({N\o r})^2\\
&\ll_{\e} M^2N^{3+2\e}\log(2M).
\end{align*}
Hence, the contribution of $R_1(T;\,h,\,k)$ is
$$
\ll{T\log T\o M^2N^{4-2\e}}\cdot M^2N^{3+2\e}\log(2M)\ll{T\log^2T\o
N^{1-8\e}}.
$$

2. The contribution of $R(T;\,h,\,k)$

By the functional equation in Lemma 14, we get
\begin{align*}
R(T;\,h,\,k)&={1\o 2\pi i}\int_{-\e-i\infty}^{-\e+i\infty}D(s;
\,{h^\ast\o k^\ast})\Gamma(s)\frak{z}^{-s}ds\\
&={1\o 2\pi i}\int_{-\e-i\infty}^{-\e+i\infty}2G^2(s)k^{\ast(1-2s)}
\Bigl(D(1-s;\,{\overline{h^\ast}\o k^\ast})\\
&\ -\cos(\pi s)D(1-s;\,-{\overline{h^\ast}\o k^\ast})\Bigr)
\Gamma(s)\frak{z}^{-s}ds\\
\quad &=k^\ast\sum_{l=1}^\infty{d(l)\o l}\cdot{1\o 2\pi
i}\int_{-\e-i\infty}
^{-\e+i\infty}2G^2(s)\cdot{l^s\o (h^\ast k^\ast)^s}\cdot\qquad\qquad
(4.14)\\
&\ \cdot\Bigl(e(l{\overline{h^\ast}\o k^\ast})-\cos(\pi s)e(-l
{\overline{h^\ast}\o k^\ast})\Bigr)\Gamma(s)\Bigl(4\pi\sin({1\o 2T})
e^{i\o 2T}\Bigr)^{-s}ds\\
&=k^\ast\sum_{l=1}^\infty{d(l)\o l}\cdot{1\o 2\pi
i}\int_{-\e-i\infty}
^{-\e+i\infty}U(s,\,T)({l\o h^\ast k^\ast})^s \cdot\\
&\ \cdot\Bigl(e(l{\overline{h^\ast}\o k^\ast})-\cos(\pi s)e(-l
{\overline{h^\ast}\o k^\ast})\Bigr)ds,
\end{align*}
where
$$
U(s,\,T)=2G^2(s)\Gamma(s)\Bigl(4\pi\sin({1\o 2T})e^{i\o 2T}\Bigr)
^{-s}. \eqno (4.15)
$$

The contribution of $R(T;\,h,\,k)$ is
$$
\ll\Bigl|\sum_{K<k,\,h\leq 8^\e K}b(k)\overline{b(h)}({h\o k})^{1\o
2}R(T;\,h,\,k)\Bigr|,
$$
while
\begin{align*}
&\ \,\sum_{K<k,\,h\leq 8^\e K}b(k)\overline{b(h)}({h\o k})^{1\o
2}R(T;\,h,\,k)\\
&=\sum_{K<k\leq 8^\e K}\sum_{K<h\leq 8^\e K}{b(k)\o k^{1\o 2}}\cdot
\overline{b(h)}h^{1\o 2}\cdot k^\ast\sum_{l=1}^\infty{d(l)\o l}\cdot\\
&\qquad \cdot{1\o 2\pi
i}\int_{-\e-i\infty}^{-\e+i\infty}U(s,\,T)({l\o h^\ast
k^\ast})^s\Bigl(e(l{\overline{h^\ast}\o k^\ast})-\cos(\pi s)e(-l{
\overline{h^\ast}\o k^\ast})\Bigr)ds\\
&=\sum_{l=1}^\infty{d(l)\o l}\cdot{1\o 2\pi i}\int_{-\e-i\infty}^
{-\e+i\infty}U(s,\,T)l^s\Bigl(\sum_{K<k\leq 8^\e K}\sum_{K<h\leq
8^\e K}{b(k)\o k^{1\o 2}}\cdot\overline{b(h)}h^{1\o 2}\cdot
(4.16)\\
&\qquad \cdot{k^\ast\o (h^\ast k^\ast)^s}e(l{\overline{h^\ast}\o
k^\ast}) -\cos(\pi s)\sum_{K<k\leq 8^\e K}\sum_{K<h\leq 8^\e
K}{b(k)\o k^{1\o 2}}\cdot\overline{b(h)}h^{1\o 2}\cdot\\
&\qquad \cdot{k^\ast\o (h^\ast k^\ast)^s}e(-l{\overline{h^\ast}\o
k^\ast})
\Bigr)ds\\
&=\sum_{l=1}^\infty{d(l)\o l}\cdot{1\o 2\pi i}\int_{-\e-i\infty}^
{-\e+i\infty}U(s,\,T)l^s(Q(l,\,s)-\cos(\pi s)Q(-l,\,s))ds,
\end{align*}
where
\begin{align*}
Q(l,\,s)&=\sum_{K<k\leq 8^\e K}\sum_{K<h\leq 8^\e K}{b(k)\o
k^{1\o 2}}\cdot\overline{b(h)}h^{1\o 2}\cdot{k^\ast\o (h^\ast
k^\ast)^s}e(l{\overline{h^\ast}\o k^\ast})\qquad\qquad (4.17)\\
&=\sum_{K<k\leq 8^\e K}\sum_{K<h\leq 8^\e K}b(k)\overline{b(h)}\cdot
{1\o (k^\ast h^\ast)^{s-{1\o 2}}}e(l{\overline{h^\ast}\o k^\ast}).
\end{align*}

For $s=-\e+it$, by the discussion in Section 5 of [4],
\begin{align*}
\qquad\qquad&U(s,\,T)l^s\ll ({T\o l})^\e (|t|+1)^{{1\o
2}+\e}\exp(({1\o 2T}-{3\o
2}\pi)|t|),\qquad\qquad\ \ (4.18)\\
&U(s,\,T)l^s\cos(\pi s)\ll ({T\o l})^\e (|t|+1)^{{1\o 2}+\e}\exp(
({1\o 2T}-{\pi\o 2})|t|).\qquad (4.19)
\end{align*}

In the following we shall estimate $Q(l,\,s)$ for $s=-\e+it$.
\begin{align*}
Q(l,\,s)&=\sum_{K<k\leq 8^\e K}\sum_{K<h\leq 8^\e
K}b(k)\overline{b(h)}\cdot{1\o (k^\ast h^\ast)^{s-{1\o 2}}}
e(l{\overline{h^\ast}\o k^\ast})\\
&=\sum_{K<k\leq 8^\e K}\sum_{K<h\leq 8^\e K}b(k)\overline{b(h)}
\cdot{(k,\,h)^{2s-1}\o (kh)^{s-{1\o 2}}}e(l{\overline{h^\ast}\o
k^\ast})\\
&=\sum_{d\leq 8^\e K}d^{2s-1}\sum_{K<k\leq 8^\e
K}\sum_{\substack{K<h\leq 8^\e K\\ (k,\,h)=d}}{b(k)\overline{b(h)}\o
(kh)^{s-{1\o 2}}}e\Bigl(l{\overline{({h\o d})}\o {k\o d}}\Bigr)
\end{align*}
\begin{align*}
&=\sum_{d\leq 8^\e K}d^{2s-1}\sum_{N<n_1\leq 2^\e N}\sum_{N<n_2\leq
2^\e N}\sum_{M<m_1\leq 2^\e M}\cdot\qquad\qquad\qquad  (4.20)\\
&\ \cdot\sum_{\substack{M<m_2\leq 2^\e M\\
(m_1n_1^2,\,m_2n_2^2)=d}} {a(n_2)\o m_2^{1\o
2}}\cdot{\overline{a(n_1)}\o m_1^{1\o 2}} \cdot{1\o
(m_1m_2n_1^2n_2^2)^{s-{1\o 2}}}e\Bigl(l{\overline{({m_1n_1^2\o
d})}\o {m_2n_2^2\o d}}\Bigr)\\
&=\sum_{d\leq 8^\e K}d^{2s-1}\sum_{N<n_1\leq 2^\e N}\sum_{N<n_2\leq
2^\e N}{a(n_2)\overline{a(n_1)}\o (n_1n_2)^{2s-1}}\cdot\\
&\ \cdot\sum_{M<m_1\leq 2^\e M}\sum_{\substack{M<m_2\leq 2^\e M\\
(m_1n_1^2, \,m_2n_2^2)=d}}{1\o
(m_1m_2)^s}e\Bigl(l{\overline{({m_1n_1^2\o d})}
\o {m_2n_2^2\o d}}\Bigr)\\
&=\sum_{d\leq 8^\e K}d^{2s-1}\sum_{N<n_1\leq 2^\e N}\sum_{N<n_2\leq
2^\e N}{a(n_2)\overline{a(n_1)}\o (n_1n_2)^{2s-1}}\cdot B(l,\,s,\,
n_1,\,n_2,\,d),
\end{align*}
where
$$
B(l,\,s,\,n_1,\,n_2,\,d)
=\sum_{\substack{M<m_2\leq 2^\e M\\ d|m_2n_2^2}}\sum_{\substack
{M<m_1\leq 2^\e M\\
(m_1n_1^2,\,m_2n_2^2)=d}}{1\o (m_1m_2)^s}
e\Bigl(l{\overline{({m_1n_1^2\o d})}\o {m_2n_2^2\o d}}\Bigr).\
(4.21)
$$

We shall estimate
$$
\sum_{\substack{M<m_1\leq M_1\\ (m_1n_1^2,\,m_2n_2^2)=d}}
e\Bigl(l{\overline{({m_1n_1^2\o d})}\o {m_2n_2^2\o d}}\Bigr)
$$
for $M<M_1\leq 2^\e M$. Let $(m_1,\,d)=d_1$. Write $d=d_1d_2$. We
see $(d_2,\,{m_1\o d_1})=1$. Hence, $d|m_1n_1^2\Longrightarrow
d_2|n_1^2\Longrightarrow d_2\leq 4^\e N^2$. By Lemma 17,
\begin{align*}
&\ \,\sum_{\substack{M<m_1\leq M_1\\ (m_1n_1^2,\,m_2n_2^2)=d}}
e\Bigl(l{\overline{({m_1n_1^2\o d})}\o {m_2n_2^2\o d}}\Bigr)\\
&=\sum_{d_1|d}\Bigl(\sum_{\substack{M<m_1\leq M_1\\ (m_1,\,d)=d_1\\
({m_1\o d_1}\cdot{n_1^2\o d_2},\,{m_2n_2^2\o d})=1}}
e\Bigl(l{\overline{({m_1n_1^2\o d})}\o {m_2n_2^2\o d}}\Bigr)\Bigr)\\
&=\sum_{d_1|d}\Bigl(\sum_{\substack{{M\o d_1}<m_1'\leq {M_1\o d_1}\\
(m_1', \,d_2)=1\\ (m_1',\,{m_2n_2^2\o d})=1\\ ({n_1^2\o
d_2},\,{m_2n_2^2\o d})=1}} e\Bigl(l{\overline{({n_1^2\o
d_2})}\cdot\overline{m_1'}\o {m_2n_2^2\o d}}\Bigr)\Bigr)
\end{align*}
\begin{align*}
&=\sum_{\substack{d_1|d\\ ({n_1^2\o d_2},\,{m_2n_2^2\o d})=1}}
\Bigl(\sum_{\substack{{M\o d_1}<m_1'\leq {M_1\o d_1}\\
(m_1',\,{m_2 n_2^2\o d})=1\\
(m_1',\,d_2)=1}}e\Bigl(l\overline{({n_1^2\o d_2})}
\cdot{\overline{m_1'}\o {m_2n_2^2\o d}}\Bigr)\Bigr)\\
&\ll\sum_{\substack{d_1|d\\ ({n_1^2\o d_2},\,{m_2n_2^2\o
d})=1}}\Bigl(l \overline{({n_1^2\o d_2})},\,{m_2n_2^2\o
d}\Bigr)^{1\o 2}\Bigl({m_2n_2^2\o d}\Bigr)^{{1\o 2}+\e}d_2^\e \\
&\ll\Bigl(l,\,{m_2n_2^2\o d}\Bigr)^{1\o
2}\sum_{d_1|d}\Bigl({m_2n_2^2 \o d}\Bigr)^{{1\o 2}+\e}d_2^\e\\
&\ll\Bigl(\sum_{d_1|d}1\Bigr)\Bigl(l,\,{m_2n_2^2\o d}\Bigr)^{1\o 2}
\Bigl({MN^2\o d}\Bigr)^{{1\o 2}+\e}d^\e \\
&\ll\Bigl(l,\,{m_2n_2^2\o d}\Bigr)^{1\o 2}\Bigl({MN^2\o d}\Bigr)^
{{1\o 2}+\e}d^{2\e},
\end{align*}
here we note $d_2\leq 4^\e N^2\Longrightarrow{M_1\o d_1}<{2m_2n_2^2
\o d}$.

By the above estimate and partial summation, for $s=-\e+it$, we have
$$
\sum_{\substack{M<m_1\leq 2^\e M\\ (m_1n_1^2,\,m_2n_2^2)=d}}{1\o
m_1^s}e\Bigl(l{\overline{({m_1n_1^2\o d})}\o {m_2n_2^2\o d}} \Bigr)
\ll(|t|+1)M^\e\Bigl(l,\,{m_2n_2^2\o d}\Bigr)^{1\o 2}\Bigl({MN^2\o
d}\Bigr)^{{1\o 2}+\e}d^{2\e}.
$$
By Lemma 16,
\begin{align*}
B(l,\,s,\,n_1,\,n_2,\,d)&\ll (|t|+1)\Bigl({MN^2\o d}\Bigr)^{{1\o 2}
+\e}M^{2\e}d^{2\e}\sum_{\substack{M<m_2\leq 2^\e
M\\ d|m_2n_2^2}}\Bigl(l,\,{m_2n_2^2\o d}\Bigr)^{1\o 2}\\
&\leq(|t|+1)\Bigl({MN^2\o d}\Bigr)^{{1\o
2}+\e}M^{2\e}d^{2\e}\sum_{M<
m_2\leq 2^\e M}(l,\,m_2)^{1\o 2}(l,\,n_2^2)^{1\o 2}\\
&\ll(|t|+1)\Bigl({MN^2\o d}\Bigr)^{{1\o
2}+\e}M^{1+2\e}d^{2\e}(l,\,n_2 ^2)^{1\o 2}l^{\e\o 4}.
\end{align*}
By Lemma 16 again, we get
\begin{align*}
Q(l,\,s)&\ll (|t|+1)(MN^2)^{{1\o 2}+\e}M^{1+2\e}\sum_{d\leq 8^\e K}
{1\o d^{{3\o 2}+\e}}\cdot\\
&\ \cdot\sum_{N<n_1\leq 2^\e N}\sum_{N< n_2\leq 2^\e
N}|a(n_1)a(n_2)|
N^{2(1+2\e)}(l,\,n_2^2)^{1\o 2}l^{\e\o 4}\\
&\ll(|t|+1)M^{{3\o 2}+3\e}N^{2+8\e}\sum_{N<n_2\leq 2^\e N}(n_2,\,l)
l^{\e\o 4}\\
&\ll (|t|+1)M^{{3\o 2}+3\e}N^{3+8\e}l^{\e\o 2}.
\end{align*}

The contribution of $R(T;\,h,\,k)$ is
\begin{align*}
&\ll\sum_{l=1}^\infty{d(l)\o l}({T\o l})^\e
\int_{-\infty}^\infty (|t|+1)^{{3\o 2}+\e}\exp(({1\o 2T}-{\pi\o
2})|t|)dt\cdot M^{{3\o 2}
+3\e}N^{3+8\e}l^{\e\o 2}\\
&\ll T^\e M^{{3\o 2}+3\e}N^{3+8\e}\sum_{l=1}^\infty{d(l)\o l^{1+{\e
\o 2}}}\\
&\ll T^\e M^{{3\o 2}+3\e}N^{3+8\e}\\
&\ll{T\o N^{1-8\e}}.
\end{align*}

3. The contribution of $R_0(T;\,h,\,k)$

Using Lemma 14, (4.12) and the estimates in 2., we get that the
contribution of $R_0(T;\,h,\,k)$ is
$$
\ll\Bigl|\sum_{K<k,\,h\leq 8^\e K}b(k)\overline{b(h)}({h\o k})^{1\o
2}R_0(T;\,h,\,k)\Bigr|,
$$
while
\begin{align*}
&\ \,\sum_{K<k,\,h\leq 8^\e K}b(k)\overline{b(h)}({h\o k})^{1\o
2}R_0(T;\,h,\,k)\\
&=\sum_{K<k\leq 8^\e K}\sum_{K<h\leq 8^\e K}b(k)\overline{b(h)}({h\o
k})^{1\o 2}\Bigl({1\o 4}-{1\o \pi
i}\sum_{a=1}^{k^\ast}\beta(a,\,k^\ast)\cdot\\
&\qquad\cdot\sum_{0<b<{k^\ast\o 2}}\eta(b,\,k^\ast)e(ab
{\overline{h^\ast}\o k^\ast})\Bigr)\\
&=\sum_{d\leq 8^\e K}\sum_{K<k\leq 8^\e K}\sum_{\substack{K<h\leq
8^\e K\\ (k,\,h)=d}}b(k)\overline{b(h)}({h\o k})^{1\o 2}\Bigl({1\o
4}-{1\o \pi i}\sum_{a=1}^{k\o d}\beta(a,\,{k\o d})\cdot\\
&\qquad\cdot\sum_{0<b<{k\o 2d}}\eta(b,\,{k\o d})e\Bigl(ab
{\overline{({h\o d})}\o {k\o d}}\Bigr)\Bigr)\\
&=\sum_{d\leq 8^\e K}\sum_{N<n_1\leq 2^\e N}\sum_{N<n_2\leq 2^\e N}
\sum_{M<m_2\leq 2^\e M}\sum_{\substack{M<m_1\leq 2^\e M\\
(m_1n_1^2,\,m_2n_2^2)=d}}{a(n_2)\o m_2^{1\o
2}}\cdot{\overline{a(n_1)}\o m_1^{1\o 2}}\cdot\\
&\ \cdot\Bigl({m_1n_1^2\o m_2n_2^2}\Bigr)^{1\o 2}\Bigl({1\o 4}-{1\o
\pi i}\sum_{a=1}^{m_2n_2^2\o d}\beta(a,\,{m_2n_2^2\o
d})\sum_{0<b<{m_2n_2^2\o 2d}}\eta(b,\,{m_2n_2^2\o
d})e\Bigl(ab{\overline{({m_1n_1^2\o d})}\o {m_2n_2^2\o
d}}\Bigr)\Bigr)\\
&=\sum_{d\leq 8^\e K}\sum_{N<n_1\leq 2^\e N}\sum_{N<n_2\leq 2^\e
N}a(n_2)\overline{a(n_1)}({n_1\o n_2})\sum_{M<m_2\leq 2^\e M}{1\o
m_2}\cdot
\end{align*}
\begin{align*}
&\ \cdot\sum_{\substack{M<m_1\leq 2^\e M\\
(m_1n_1^2,\,m_2n_2^2)=d}}\Bigl({1\o 4}-{1\o \pi
i}\sum_{a=1}^{m_2n_2^2\o d}\beta(a,\,{m_2n_2^2\o
d})\sum_{0<b<{m_2n_2^2\o 2d}}\eta(b,\,{m_2n_2^2\o
d})e\Bigl(ab{\overline{({m_1n_1^2\o d})}\o {m_2n_2^2\o
d}}\Bigr)\Bigr)\\
 &\ll\sum_{d\leq 8^\e K}\sum_{N<n_1\leq 2^\e
N}\sum_{N<n_2\leq 2^\e N}|a(n_1)a(n_2)|\sum_{M<m_2\leq 2^\e M}{1\o
m_2}\sum_{\substack{M<m_1 \leq 2^\e M\\ (m_1n_1^2,\,m_2n_2^2)=d}}1\\
&\ +\sum_{d\leq 8^\e K}\sum_{N<n_1\leq 2^\e N}\sum_{N<n_2\leq 2^\e
N}|a(n_1)a(n_2)|\sum_{\substack{M<m_2\leq 2^\e M\\ d|m_2n_2^2}}{1\o
m_2}\cdot\\
&\ \cdot\sum_{a=1}^{m_2n_2^2\o d}|\beta(a,\,{m_2n_2^2\o
d})|\sum_{0<b<{m_2n_2^2\o 2d}}\eta(b,\,{m_2n_2^2\o
d})\Bigl|\sum_{\substack {M<m_1\leq 2^\e M\\
(m_1n_1^2,\,m_2n_2^2)=d}}e\Bigl(ab{\overline{({m_1n_1^2\o d})}\o
{m_2n_2^2\o d}}\Bigr)\Bigr|\\
&\ll\sum_{N<n_1\leq 2^\e N}\sum_{N<n_2\leq 2^\e N}\sum_{M<m_1\leq
2^\e M}\sum_{M<m_2 \leq 2^\e M}{|a(n_1)a(n_2)|\o m_2}\\
&\ +\sum_{d\leq 8^\e K}\sum_{N<n_1\leq 2^\e N}\sum_{N<n_2\leq 2^\e
N}|a(n_1)a(n_2)|\sum_{\substack{M<m_2\leq 2^\e M\\ d|m_2n_2^2}}{1\o
m_2}\cdot\\
&\ \cdot\sum_{a=1}^{m_2n_2^2\o d}|\beta(a,\,{m_2n_2^2\o
d})|\sum_{0<b<{m_2n_2^2\o 2d}}{1\o b}\Bigl(ab,\,{m_2n_2^2\o
d}\Bigr)^{1\o 2}\Bigl({MN^2\o d}\Bigr)^{{1\o 2}+\e}d^{2\e}\\
&\ll N^{-2+2\e}M^{-1}(MN)^2+N^{-2+2\e}M^{-1}\sum_{d\leq 8^\e
K}\sum_{N<n_1\leq 2^\e N}\sum_{N<n_2\leq 2^\e N}\cdot\\
&\ \cdot\sum_{\substack{M<m_2\leq 2^\e M\\ d|m_2n_2^2}}\Bigl({MN^2\o
d}\Bigr)^{{1\o 2}+\e}d^{2\e}\sum_{a=1}^{m_2n_2^2\o
d}|\beta(a,\,{m_2n_2^2\o d})|\Bigl(a,\,{m_2n_2^2\o d}\Bigr)^{1\o
2}\cdot\\
&\ \cdot\sum_{0<b<{m_2n_2^2\o 2d}}{1\o b}\Bigl(b,\,{m_2n_2^2\o
d}\Bigr)^{1\o 2}.
\end{align*}

We have
$$
\sum_{0<b<{m_2n_2^2\o 2d}}{1\o b}\Bigl(b,\,{m_2n_2^2\o d}\Bigr)^{1\o
2}=\sum_{r|{m_2n_2^2\o d}}r^{1\o 2}\sum_{\substack{0<b<{m_2n_2^2\o
2d}\\ (b,\,{m_2n_2^2\o d})=r}}{1\o b}
$$
\begin{align*}
&\leq\sum_{r|{m_2n_2^2\o d}}r^{1\o 2}\sum_{\substack{0<b\leq
{m_2n_2^2\o 2d}\\ r|b}}{1\o b}\\
&\ll\sum_{r|{m_2n_2^2\o d}}{1\o r^{1\o 2}}\log({2m_2n_2^2\o d})\\
&\ll \Bigl({m_2n_2^2\o d}\Bigr)^{\e\o 4}\\
&\ll \Bigl({MN^2\o d}\Bigr)^{\e\o 4}
\end{align*}
and
\begin{align*}
&\ \,\sum_{a=1}^{m_2n_2^2\o d}|\beta(a,\,{m_2n_2^2\o
d})|\Bigl(a,\,{m_2n_2^2\o d}\Bigr)^{1\o 2}\\
&\ll{m_2n_2^2\o d}+\sum_{1\leq a\leq{m_2n_2^2\o 2d}}{{m_2n_2^2\o
d}\o a}\Bigl(a,\,{m_2n_2^2\o d}\Bigr)^{1\o 2}\\
&\ +\sum_{{m_2n_2^2\o 2d}<a\leq{m_2n_2^2\o d}-1}{{m_2n_2^2\o d}\o
{m_2n_2^2\o d}-a}\Bigl({m_2n_2^2\o d}-a,\,{m_2n_2^2\o d}\Bigr)^{1\o 2}\\
&\ll{m_2n_2^2\o d}+{m_2n_2^2\o d}\sum_{1\leq a\leq{m_2n_2^2\o 2d}}{1\o
a}\Bigl(a,\,{m_2n_2^2\o d}\Bigr)^{1\o 2}\\
&\ll\Bigl({MN^2\o d}\Bigr)^{1+{\e\o 4}}.
\end{align*}

Therefore the contribution of $R_0(T;\,h,\,k)$ is
\begin{align*}
&\ll MN^{2\e}+N^{-2+2\e}M^{-1}\sum_{d\leq 8^\e K}\sum_{N<n_1\leq
2^\e N}\sum_{N<n_2\leq 2^\e N}\cdot\\
&\ \cdot\sum_{M<m_2\leq 2^\e M}\Bigl({MN^2\o d}\Bigr)^{{1\o
2}+\e}d^{2\e}
\Bigl({MN^2\o d}\Bigr)^{1+{\e\o 2}}\\
&\ll MN^{2\e}+N^{2\e}(MN^2)^{{3\o 2}+{3\o 2}\e}\sum_{d\leq 8^\e
K}{1\o d^{{3\o 2}-{\e\o 2}}}\\
&\ll MN^{2\e}+M^{{3\o 2}+{3\o 2}\e}N^{3+5\e}\\
&\ll {T\o N^{1-8\e}}.
\end{align*}

Combining all of the above, we get
$$
\int_0^\infty e^{-{t\o T}}|\zeta({1\o 2}+it)|^2\Bigl|\sum_{K< k\leq
8^\e K}{b(k)\o k^{it}}\Bigr|^2\ll{T\log^2T\o N^{1-8\e}}.
$$
Hence,
\begin{align*}
&\log T\sum_{r\leq{1\o 2\e\log 2}\log T+O(1)}\int_0^\infty
e^{-{t\o T}}|\zeta({1\o 2}+it)|^2\cdot\\
&\cdot\Bigl|\sum_{2^{\e r}< m\leq 2^\e\cdot 2^{\e r}}{1\o m^{{1\o
2}+it}}\Bigr|^2\Bigl|\sum_{N< n\leq 2^\e N}{a(n)\o
n^{2it}}\Bigr|^2dt \ll{T\log^4T\o N^{1-8\e}}.
\end{align*}

So far the proof of Proposition 3 is finished.

Proof of Proposition 2. We observe that the measure of the set of
all $t$ such that $\frac {T}{2} \leq t \leq T$ and $2t \in U_4$ is
$\ll T^{11\e}(\log T)^{10}$. We suppose firstly that $k=2m$ with
positive integer $m$. By Proposition 1, Lemmas 8, 9 and 11,
\begin{align*}
&\ \,\int_{T\o 2}^T|\zeta({1\o 2}+it)|^4|\zeta(1+2it)|^{-2m}dt\\
&=\int_{{T\o 2}\leq t\leq T,\,2t\not\in U_4}|\zeta({1\o 2}+it)|^4|
\zeta(1+2it)|^{-2m}dt\\
&\ +\int_{{T\o 2}\leq t\leq T,\,2t\in U_4}|\zeta({1\o 2}+it)|^4
|\zeta(1+2it)|^{-2m}dt\\
&\ll\int_{{T\o 2}\leq t\leq T,\,2t\not\in U_4}|\zeta({1\o 2}+it)|^4
\Bigl(\Bigl|\sum_{l\leq X}{\mu(l)\o l^{1+2it}}\cdot e^{-{l\o X}}
\Bigr|^{2m}+ O(1)\Bigr)dt+O(T\log^4T)\\
&\ll\int_{T\o 2}^T|\zeta({1\o 2}+it)|^4\Bigl|\sum_{l\leq X}{\mu(l)\o
l^{1+2it}}\cdot e^{-{l\o X}}\Bigr|^{2m}dt+O(T\log^4T)\\
&=\int_{T\o 2}^T|\zeta({1\o 2}+it)|^4\Bigl|\sum_{n\leq X^m}{a(n)\o
n^{2it}}\Bigr|^2 dt+O(T\log^4T),
\end{align*}
where
$$
a(n)={1\o n}\sum_{l_1\cdots l_m=n}\mu(l_1)\cdots\mu(l_m)\exp(-{(l_1
+\cdots+l_m)\o X}),
$$
$U_4$ is defined as in (3.2), $X$ is defined as in (3.8). We can see
$$
X^m=\exp({2m\o \e}(\log T)^{3\o 4})\ll T^{{1\o 16}-2\e}
$$
and
$$
a(n)=O(n^{-1+\e}).
$$
By Cauchy's inequality,
\begin{align*}
\Bigl|\sum_{n\leq X^m}{a(n)\o n^{2it}}\Bigr|^2
&=\Bigl|\sum_{s\leq{{m\log X\o \e\log 2}-1}}{1\o 2^{\e s\o 4}}\cdot
2^{\e s\o 4}\sum_{2^{\e s}<n\leq 2^\e 2^{\e s}}{a(n)\o n^{2it}}
\Bigr|^2\\
&\leq\sum_{s\leq{{m\log X\o \e\log 2}-1}}{1\o 2^{\e s\o 2}}
\sum_{s\leq{{m\log X\o \e\log 2}-1}}2^{\e s\o 2}\Bigl|\sum_{2^{\e
s}<n\leq 2^\e 2^{\e s}}{a(n)\o n^{2it}}\Bigr|^2\\
&\ll\sum_{s\leq{{m\log X\o \e\log 2}-1}}2^{\e s\o 2}\Bigl|\sum_{2^
{\e s}<n\leq 2^\e 2^{\e s}}{a(n)\o n^{2it}}\Bigr|^2.
\end{align*}
Hence, Proposition 3 yields
\begin{align*}
&\ \,\int_{T\o 2}^T|\zeta({1\o 2}+it)|^4\Bigl|\sum_{n\leq
X^m}{a(n)
\o n^{2it}}\Bigr|^2dt\\
&\ll\sum_{s\leq{{m\log X\o \e\log 2}-1}}2^{\e s\o 2}\int_{T\o 2}^T
|\zeta({1\o 2}+it)|^4\Bigl|\sum_{2^{\e s}<n\leq 2^\e 2^{\e s}}{a(n)
\o n^{2it}}\Bigr|^2\\
&\ll\sum_{s\leq{{m\log X\o \e\log 2}-1}}2^{\e s\o 2}\cdot{T\log^4T
\o 2^{\e s(1-8\e)}}\\
&\ll T\log^4T.
\end{align*}
Thus,
$$
\int_{T\o 2}^T|\zeta({1\o 2}+it)|^4|\zeta(1+2it)|^{-2m}dt\ll
T\log^4T.
$$
Therefore
$$
\int_1^T|\zeta({1\o 2}+it)|^4|\zeta(1+2it)|^{-2m}dt\ll T\log^4T.
$$

For the general $k>0$, we have an even integer $2m$ such that
$k<2m$. By H\"older's inequality,
\begin{align*}
&\ \,\int_1^T|\zeta({1\o 2}+it)|^4|\zeta(1+2it)|^{-k}dt\\
&=\int_1^T|\zeta({1\o 2}+it)|^{4\cdot{2m-k\o 2m}}\cdot|\zeta({1\o 2}
+it)|^{4\cdot{k\o 2m}}|\zeta(1+2it)|^{-k}dt\\
&\leq\Bigl(\int_1^T|\zeta({1\o 2}+it)|^{4\cdot{2m-k\o 2m}\cdot{2m\o
2m-k}}dt\Bigr)^{2m-k\o 2m}\cdot\\
&\ \cdot\Bigl(\int_1^T|\zeta({1\o 2}+it)|^{4\cdot{k\o 2m}\cdot{2m \o
k}}|\zeta(1+2it)|^{-k\cdot{2m\o k}}dt\Bigr)^{k\o 2m}
\end{align*}
\begin{align*}
&=\Bigl(\int_1^T|\zeta({1\o 2}+it)|^4dt\Bigr)^{2m-k\o
2m}\Bigl(\int_1
^T|\zeta({1\o 2}+it)|^4|\zeta(1+2it)|^{-2m}dt\Bigr)^{k\o 2m}\\
&\ll T\log^4T.
\end{align*}

So far the proof of Proposition 2 is finished.

\vskip.3in
\noindent{\bf 5. The proof of Theorem}

We shall apply Lemma 5. For ${\rm Re}(s)>1$, let
$$
f(s)=\sum_{n=1}^\infty{d^2(n)\o n^s}.
$$
By Lemma 6,
$$
f(s)={\zeta^4(s)\o \zeta(2s)}.
$$
It is easy to see that $\psi(n)=n^\e$ which is non-decreasing. As
$\sigma\rightarrow 1^+$,
$$
\sum_{n=1}^\infty{d^2(n)\o n^\sigma}={\zeta^4(\sigma)\o \zeta(2
\sigma)}=O({1\o (\sigma-1)^4}).
$$
Let $c=1+\e,\,Y=[x]+{1\o 2},\,T=x^{3\o 4}$. Then
\begin{align*}
\sum_{n\leq x}d^2(n)&=\sum_{n<Y}d^2(n)+O_{\e}(x^\e)\\
&={1\o 2\pi i}\int_{1+\e-iT}^{1+\e+iT}{\zeta^4(s)\o \zeta(2s)}
\cdot{Y^s\o s}ds+O_{\e}(x^{{1\o 4}+2\e}).
\end{align*}
We move the line of integration to ${\rm Re}(s)={1\o 2}$. The residue
of ${\zeta^4(s)\o\zeta(2s)} \cdot{Y^s\o s}$ at $s=1$ is
$$
YP(\log Y)=xP(\log x)+O(x^\e).
$$
By Lemmas 8 and 10,
$$
{1\o 2\pi i}\int_{{1\o 2}+iT}^{1+\e+iT}{\zeta^4(s)\o \zeta(2s)}
\cdot{Y^s\o s}ds \ll\max_{{1\o 2}\leq\sigma\leq 1+\e}T^{{4\o
3}(1-\sigma)+4\e} \log T\cdot{x^\sigma\o T}\ll x^{{1\o 4}+4\e}.
$$
In the same way,
$$
{1\o 2\pi i}\int_{{1\o 2}-iT}^{1+\e-iT}{\zeta^4(s)\o \zeta(2s)}
\cdot{Y^s\o s}ds\ll x^{{1\o 4}+4\e}.
$$
Hence,
$$
E(x)={1\o 2\pi i}\int_{{1\o 2}-iT}^{{1\o 2}+iT}{\zeta^4(s)\o
\zeta(2s)}\cdot{Y^s\o s}ds+O(x^{{1\o 4}+4\e}).
$$
It follows from Proposition 2 that,
\begin{align*}
E(x)&\ll x^{1\o 2}\sum_{k\leq{\log T\o \log 2}}\int_{2^{k-1}}^{2^k}
{|\zeta({1\o 2}+it)|^4\o |\zeta(1+2it)|}{dt\o t}+O(x^{1\o 2})\\
&\ll x^{1\o 2}\sum_{k\leq{\log T\o \log 2}}{1\o 2^k}\int_1^{2^k}
{|\zeta({1\o 2}+it)|^4\o |\zeta(1+2it)|}dt+O(x^{1\o 2})\\
&\ll x^{1\o 2}\sum_{k\leq{\log T\o \log 2}}{1\o 2^k}\cdot 2^kk^4\\
&\ll x^{1\o 2}\log^5x.
\end{align*}

Thus, the proof of the Theorem is complete.

\vskip.3in
\noindent{\bf 6. Some remarks}

By the method of this paper, we can prove that if $k$ is any given
positive number, $a$ is a given positive integer, then
$$
\int_1^T{|\zeta({1\o 2}+it)|^4\o |\zeta(1+ait)|^k}dt\ll_{k,\,a}
T\log^4T.
$$

We note that if ${\rm Re}(s)>1$,
$$
\sum_{n=1}^\infty{d(n^3)\o n^s}={\zeta^4(s)\o \zeta^3(2s)}\cdot
G_1(s), \eqno (6.1)
$$
where
$$
G_1(s)=\prod_p{(1+{2\o p^s})\o (1-{1\o p^s})(1+{1\o p^s})^3},
$$
$G_1(s)$ is absolutely convergent for ${\rm Re}(s)>{1\o 3}$. One can
see page 95 in [2]. Using the method similar to that in this paper,
we can prove the following proposition.

{\bf Proposition 4}. We have
$$
\sum_{n\leq x}d(n^3)=xP_1(\log x)+O(x^{1\o 2}(\log x)^5),
$$
where $P_1(y)$ is a suitable cubic polynomial in $y$.

In 2006, M. Z. Garaev, F. Luca and W. G. Nowak[2] proved that as
$x\rightarrow \infty$, if $y=y(x)$ satisfies
$$
{y\o x^{1\o 2}\log x\log\log x}\rightarrow\infty,
$$
then
$$
\sum_{x<n\leq x+y}d^2(n)\sim{6\o \pi^2}y(\log x)^3,
$$
and that as $x\rightarrow\infty$, if $y=y(x)$ satisfies
$$
{y\o x^{1\o 2}\log x(\log\log x)^3}\rightarrow\infty,
$$
then
$$
\sum_{x<n\leq x+y}d(n^3)\sim B_0y(\log x)^3,
$$
where $B_0$ is a positive constant.

Combining the method of this paper with that of [2], we can prove
the following proposition.

{\bf Proposition 5}. As $x\rightarrow\infty$, if $y=y(x)$ satisfies
$$
{y\o x^{1\o 2}\log x}\rightarrow\infty,
$$
then
$$
\sum_{x<n\leq x+y}d^2(n)\sim{6\o \pi^2}y(\log x)^3
$$
and
$$
\sum_{x<n\leq x+y}d(n^3)\sim B_0y(\log x)^3.
$$

Let $r(n)$ be the number of representations of $n$ as the sum of two
squares. In 2004, M. K\"uhleitner and W. G. Nowak[6] proved that
$$
\sum_{n\leq x}r^2(n)=4x\log x+B_1x+O(x^{1\o 2}(\log x)^3\log\log x),
\eqno (6.2)
$$
where $B_1$ is a positive constant, and that
$$
\sum_{n\leq x}r(n^3)=A_2x\log x+B_2x+O(x^{1\o 2}(\log x)^3(\log \log
x)^2), \eqno (6.3)
$$
where $A_2,\,B_2$ are positive constants.

Let ${\bf K}={\Bbb Q}(i)$, $\zeta_{\bf K}(s)$ be the Dedekind
$\zeta$ function in the field ${\bf K}$. If ${\rm Re}(s)>1$,
\begin{align*}
\sum_{n=1}^\infty{r^2(n)\o n^s}&={16\zeta_{\bf K}^2(s)\o (1+2^{-s})
\zeta(2s)},\\
\sum_{n=1}^\infty{r(n^3)\o n^s}&={\zeta_{\bf K}^2(s)\o \zeta(2s)
\zeta_{\bf K}(2s)}\cdot G_2(s),
\end{align*}
where $G_2(s)$ is holomorphic and bounded for ${\rm Re}(s)>{1\o
3}+\e$. One can see (4.1) and (4.4) in [6].

If the result of Iwaniec[4] could be generalized to $\zeta_{\bf
K}(s)$, then the error terms in (6.2) and (6.3) could be improved to
$O(x^{1\o 2}(\log x)^3)$. Furthermore, the sums studied in [2]
$$
\sum_{x<n\leq x+y}r^2(n),\qquad\sum_{x<n\leq x+y}r(n^3),\qquad
\sum_{x<n\leq x+y}d(n)r(n)
$$
could also be improved correspondingly.

\vskip.3in \noindent{\bf Acknowledgements}

The authors would like to thank the referee for his/her nice
comments and suggestions. Chaohua Jia is supported by the National
Natural Science Foundation of China (Grant No. 11371344) and the
National Key Basic Research Program of China (Project No.
2013CB834202).


\bigskip

\

Chaohua Jia

Institute of  Mathematics, Academia Sinica, Beijing 100190, P. R.
China

Hua Loo-Keng Key Laboratory of Mathematics, Chinese Academy of
Sciences, Beijing 100190, P. R. China

E-mail: jiach@math.ac.cn

\

Ayyadurai Sankaranarayanan

School of Mathematics, Tata Institute of Fundamental Research, Homi
Bhabha Road, Mumbai 400005, India

E-mail: sank@math.tifr.res.in
\end{document}